\documentclass[11pt]{article}

\usepackage{times}
\usepackage{graphicx}
\usepackage{psfrag}
\usepackage{amsmath,amsthm}
\usepackage{amsfonts}

\theoremstyle{plain} 

\theoremstyle{definition}

\theoremstyle{remark}

\DeclareMathOperator{\Tr}{Tr}

\DeclareMathOperator{\rank}{rank}

\begin{document}

\parskip=10pt

\flushbottom 

\title{A network that learns Strassen multiplication} 

\author{
Veit Elser\\
Cornell University\\
} 

\date{}

\maketitle

\begin{abstract}
We study neural networks whose only non-linear components are multipliers, to test a new training rule in a context where the precise representation of data is paramount. These networks are challenged to discover the rules of matrix multiplication, given many examples. By limiting the number of multipliers, the network is forced to discover the Strassen multiplication rules. This is the mathematical equivalent of finding low rank decompositions of the $n\times n$ matrix multiplication tensor, $M_n$. We train these networks with the conservative learning rule, which makes minimal changes to the weights so as to give the correct output for each input at the time the input-output pair is received. Conservative learning needs a few thousand examples to find the rank 7 decomposition of $M_2$, and $10^5$ for the rank 23 decomposition of $M_3$ (the lowest known). High precision is critical, especially for  $M_3$, to discriminate between true decompositions and ``border approximations".

\end{abstract}

\section{Introduction}

The current surge of activity in machine learning with neural networks, far from being a fad, is a response to hard empirical evidence. The powers of discrimination, generalization, and abstraction exhibited by these networks matches and often outperforms humans. At the same time, from a theoretical perspective, we know that the core representation of knowledge by these networks is only approximate. The representation of images as compositions of non-negative features, for example, works well only until the images become unnatural. While synthetic data avoids this shortcoming, the allure of the network is diminished when the learned knowledge is not new or surprising.

Questions related to the Strassen algorithm \cite{S} for matrix multiplication provide an opportunity for testing neural networks in the interesting setting where the representation of knowledge is both exact and poorly understood by humans (an open mathematical problem). 
In this short contribution we define Strassen multiplication (SM) in engineering terms, construct a neural network that has the capacity to implement SM, introduce a simple protocol for training the network called conservative learning, and finish with the mandatory presentation of spectacular results.

\section{What is Strassen multiplication?}

In 1969 Volker Strassen published a paper \cite{S} that showed, among other things, that a pair of $2\times 2$ matrices could be multiplied with just seven scalar multiplications, one less than the eight required in the standard algorithm. Strassen's trick becomes practical on very large matrices where it can be applied recursively on blocks, and finding the fewest number of scalar multiplications for general $n\times n$ matrices turns out to be a deep mathematical problem. A glance at the network that implements Strassen multiplication (SM) in Figure 1 quickly establishes two things that make this application of machine learning unique:
\begin{itemize}
\item The architecture of the network and the operations performed by its components are \textit{a perfect match to the mathematical problem}, as opposed to being just a well motivated, robust platform for general machine learning applications.
\item The design of the network is deliberately naive and lives up to the standard that the machine should have the capacity \textit{to learn something that we don't already know}. Our network, in fact, does not even know how numbers are arranged to form a matrix!
\end{itemize}

\begin{figure}[t]
\begin{center}
\includegraphics[width=5in]{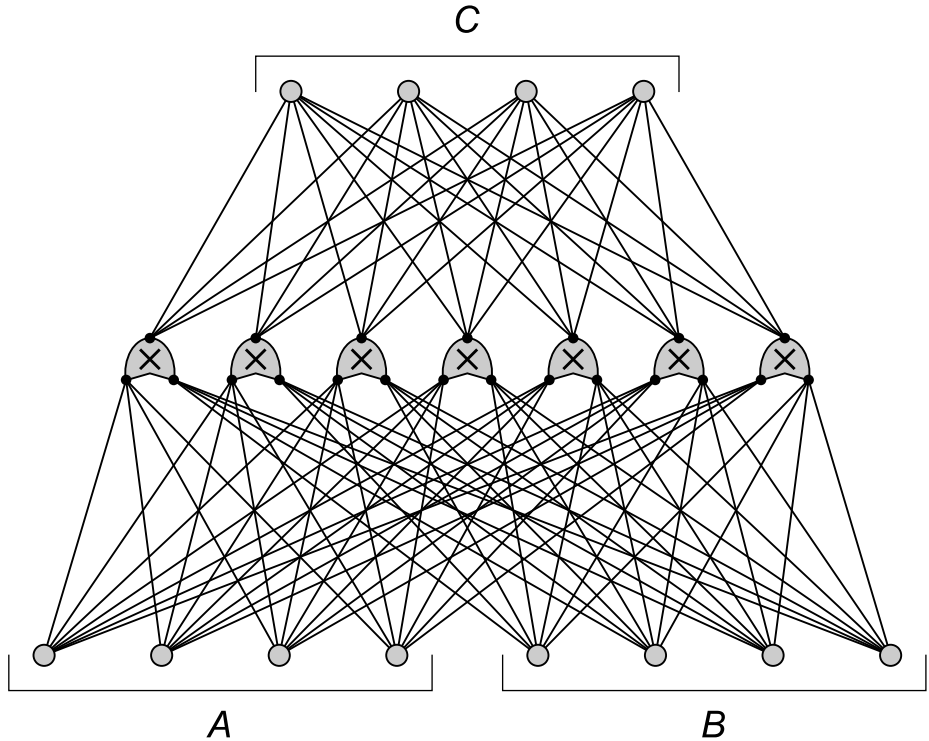}
\end{center}
\caption{A network for multiplying $2\times 2$ matrices that uses only seven multipliers. Multiplication by a constant (``weight") occurs at each of the lines connecting input and output registers (for the matrices $A$, $B$ and $C$) to the multipliers in the middle layer of the network.}
\end{figure}

Figure 1 is an engineer's representation of SM. The elements of the $2\times 2$ matrices $A$ and $B$ are read in from two registers of numbers $a(i)$ and  $b(i)$, $i=1,\ldots, 4$ on the left and right. Those from $A$ are `pooled' into the left inputs of the seven multipliers in the middle layer of the network, those from $B$ into the right inputs. Pooling is the most general linear-homogeneous operation. Denoting the inputs to the multipliers $a^*(j)$ and  $b^*(j)$, $j=1,\ldots, 7$, the pooling equations are
\begin{equation}\label{wawb}
a^*(j)=\sum_{i=1}^4 W_a(j,i)a(i)\qquad b^*(j)=\sum_{i=1}^4 W_b(j, i)b(i).
\end{equation}
The $7\times 4$ sets of real numbers $W_a$ and $W_b$ are called \textit{weights}, though they are not required to be non-negative. Multiplication by a weight is not counted as one of the multiplies in SM. Unlike the multiplications performed by the seven multipliers, weight multiplication is a linear operation in the inputs to the network. The challenge of SM is to find a set of weights that work for all conceivable matrices we wish to multiply. We will see that this set is not unique, and there exist transformations from one set to another that works just as well. In the most practical set found by Strassen \cite{S}, all the nonzero weights are just $\pm 1$, so the weight-multiplications are really just additions and subtractions.

The output end of the network works in the same way as the input end: the outputs of the seven multipliers are pooled into the four output registers of the $C$ matrix, $c(i)$, $i=1,\ldots, 4$, again in the most general way with a set of weights:
\begin{equation}\label{wc}
c(i)=\sum_{j=1}^7 W_c(i,j)\,a^*(j) b^*(j).
\end{equation}
Strassen's discovery can now be expressed in purely engineering terms. That is, there exists a fixed set of weights $W_a$, $W_b$ and $W_c$ such that the four numbers in the product of any pair of matrices can be computed with the 7-multiplier network of Figure 1 (\textit{i.e.} the 8-multiplier network implicit in the standard algorithm is sub-optimal).

\section{Mathematical digression: tensor rank}

Combining \eqref{wawb} and \eqref{wc} we obtain a curious statement about matrix multiplication, now for $n\times n$ matrices:
\begin{equation}
c(i)=\sum_{k=1}^{n^2} \sum_{l=1}^{n^2} M_n(i,k,l)\, a(k) b(l)
\end{equation}
\begin{equation}
M_n(i,k,l)=\sum_{j=1}^r W_c(i,j)\, W_a(j,k)\, W_b(j,l).
\end{equation}
The universal three-index set of numbers $M_n$ is a tensor of order three, the $(n\times n)$-\textit{matrix multiplication tensor}. We see that $M_n$ can be written as a sum of $r$ products of three 1-tensors, the atoms of all tensors. The minimum $r$ in the decomposition into products of 1-tensors is called the rank of the tensor. Unlike the problem of determining the rank of a matrix (a 2-tensor), determining the ranks of higher order tensors such as $M_n$ is computationally difficult \cite{H}. While we know $\rank{(M_2)}=7$ \cite{W}, already for the next case we currently only have bounds \cite{B,JDL}: $19\le \rank{(M_3)}\le 23$.

A general transformation of the weights that gives another correct network (decomposition of $M_n$) can be inferred from the fact that the transformed matrices
\begin{equation}\label{trans}
\tilde{A}=UAV^{-1}\qquad \tilde{B}=VBW^{-1}\qquad \tilde{C}=UCW^{-1},
\end{equation}
for arbitrary invertible $(U, V, W)$, satisfy the matrix product property whenever $A$, $B$ and $C$ do ($\tilde{A}\tilde{B}=\tilde{C}$). But feeding $\tilde{A}$ and $\tilde{B}$ into the network, and seeing the correct output $\tilde{C}$, is equivalent to feeding in $A$ and $B$ and seeing $C$ in the output after the three sets of weights have been appropriately transformed by the linear relations \eqref{trans}.

When we test our network for ``matrix multiplication fidelity" we would like to do better than spot check its accuracy on sample instances. To assess how well it will perform on \textit{any} instance of matrix multiplication we compute $\epsilon$, where
\begin{equation}\label{epsilon}
\epsilon^2=\frac{1}{(n^2)^3}\sum_{i\, k\, l}\left(M_n(i,k,l)-\sum_{j =1}^rW_c(i,j)\, W_a(j,k)\, W_b(j,l)\right)^2
\end{equation}
is the mean-square error in how our weights are decomposing the true $M_n$.

\section{Conservative learning}\label{sec:CL}

The truly amazing thing about the network in Figure 1 is not that, given the special Strassen weights, it manages to produce the product $C=A B$ with only seven multiplies. What is really remarkable is the fact that the network can learn a correct set of weights just by being shown enough examples of correctly multiplied matrices.

A huge industry has grown up around the problem of training networks. Most methods start with completely random weights and adjust these to minimize the discrepancy between the correct output and the actual output of the network. There are many schemes for minimization, though most are local and based on gradients of some measure of discrepancy (``loss") with respect to the weights. We will use a somewhat different approach for the Strassen network called conservative learning.

Conservative learning is not the oxymoron it would seem, but a combination of two very reasonable principles:
\begin{itemize}
\item When presented with a training item for which the network gives the wrong output, change the weights so that at least this item produces the correct output.
\item Make the smallest possible change to the weights when learning each new item.
\end{itemize}
By making the smallest possible change, when obsessing over the most recent item, we stand a better chance of not corrupting the accumulated knowledge derived from all previous items. To add mathematical support to this statement we analyze a very simple network.

Consider a network that implements a linear transformation from a set of inputs $x\in \mathbb{R}^m$ to a set of outputs $y\in \mathbb{R}^n$. We use matrix notation in what follows. The transformation to be learned corresponds to a matrix of weights $W^\star$, and training samples are pairs $(x,y)$ where $y=W^\star x$. Suppose $W$ is the current matrix and $(x,W^\star x)$ a new training item. As our network is linear, we choose to train only with normalized inputs, $x^T x=1$. We wish to find the smallest change $\Delta$ such that
\begin{equation}
(W+\Delta)x=W^\star x.
\end{equation}
Taking the Frobenius norm $\|\Delta\|=\Tr \Delta^T\Delta$ to define the smallness of the change, some simple linear algebra gives us the conservative learning rule:
\begin{equation}
W\to W'=W+(W^\star x-W x) x^T.
\end{equation}
The current weights are incremented by a rank 1 matrix comprising the input vector $x$ as one factor and the output discrepancy $W^\star x-W x$ as the other. The new weights $W'$ are now correct on the item $(x,W^\star x)$.

To better understand how we are doing relative to the corpus of all possible items, we make note of the following relationship between old and new weights:
\begin{eqnarray}
W'-W^\star&=&(W-W^\star)(1-x x^T)\\
&=&P_x^\perp(W-W^\star).\label{rowproj}
\end{eqnarray}
The operator $P_x^\perp$ projects rows to the orthogonal complement of $x$. From \eqref{rowproj} we infer the norm inequality
\begin{equation}
\|W'-W^\star\|\le \|W-W^\star\|,
\end{equation}
where we get equality only in the rare case that $x$ is orthogonal to all the rows of $W-W^\star$. This inequality proves that conservatism is sound: by minimally accommodating each new item there is monotone improvement in our approximation of the transformation $(W^\star$) itself.

Conservative learning has the nice feature that there are no parameters that have to be tuned, such as the step size in gradient descent. This continues to be true when the principle is applied to general networks. However, multiple layers and nonlinear components make the problem of finding the exact weight modifications intractable. Fortunately, there is a systematic procedure for keeping track of the order of smallness of the changes in the general network so that an approximate learning rule can be written down. The update procedure for the weights in the multilayer setting has similarities with the back-propagation rules in gradient-based learning, but there are also some new modes of propagation. Not surprisingly, proving convergence is also out of reach. Our experiments with the Strassen network show spectacular convergence, even without imposing a small parameter to keep changes in check.

A derivation of the conservative learning rules for the Strassen network are given in the appendix. This network has all the features that need to be addressed in general networks and therefore serves as a good vehicle for explaining the method.

\section{Conservative learning with the Strassen network}

There are two sources of randomness when training a network: the distribution of the initial weights and the distribution of the training items. We did not explore this dimension in our experiments with the Strassen network. The weights were initialized with independent uniform samples from a symmetric interval about the origin. Clearly the scale of this interval is not arbitrary, since the tensor $M_n$ being decomposed has a definite scale. A scale for the weights is also implied by the error objective \eqref{epsilon}.  Conservative learning proved not to be very sensitive to the scale of the initial weights, with only extreme limits (small and large scales) showing performance degradation. For simplicity we therefore chose to always initialize with samples drawn from $[-1,1]$. 

Adding an interesting bias to the distribution of training matrices is also something we did not try. Our input matrices were simply constructed from uniform samples of $[-1,1]$ in each element of $A$ and $B$. After rescaling to conform with our normalization convention $\Tr A^T A=\Tr B^T B=1$, these together with $C=A B$ were handed to the network for training. 

\begin{figure}[t]
\begin{center}
\includegraphics[width=4.in]{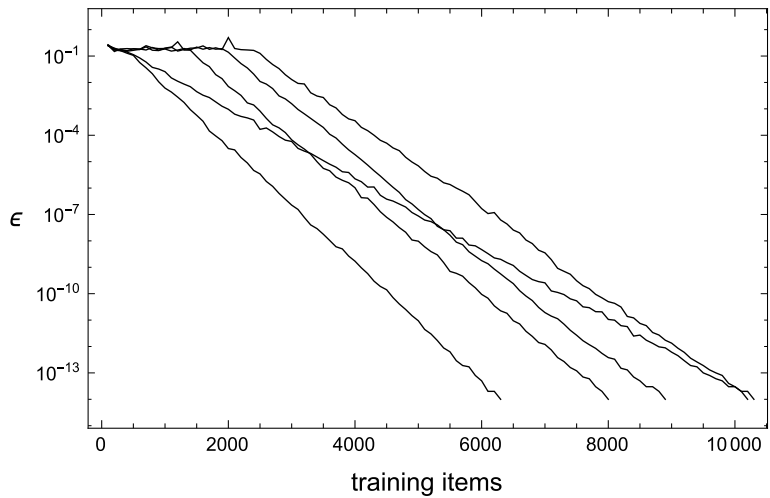}
\end{center}
\caption{Root-mean-square error $\epsilon$ in the rank 7 decomposition of the $2\times 2$ matrix multiplication tensor by the network in Figure 1. Shown are results for five runs of the network, each starting with random weights and trained on streams of random instances of correctly multiplied matrices.}
\end{figure}

It astonished us how quickly the network discovers Strassen's trick for $2\times 2$ matrices, when the weights are updated by the conservative learning rules. Figure 2 shows the decomposition error $\epsilon$ converging linearly, after an initial training set of only a few thousand. The five runs shown have slightly different convergence rates, a variability that can only be attributed to the initial weights and the early training items (before the onset of linear convergence). We did not compare with other learning rules, all of which, unlike conservative learning, come with parameters and batch protocols that we were not prepared to manage with competence.

When we remove one multiplier from the middle layer of the network in Figure 1, that is, attempt to find a rank 6 decomposition of $M_2$, the error $\epsilon$ fluctuates indefinitely. This of course is what we expect, since there is a proof that $M_2$ has rank 7 \cite{W}. For $3\times 3$ matrices the network exhibits a new kind of behavior. Since it is known that $\rank{(M_3)}\le 23$ \cite{JDL}, convergent behavior is possible in principle when we have 23 multipliers. We find that this indeed happens in about 64\% of all runs, an example of which is the lower plot in Figure 3. In the remaining runs $\epsilon$ decays much more slowly; an example of this behavior is the higher plot in Figure 3.

The phenomenon of there being two qualitatively different asymptotic states of the network, when decomposing $M_3$, is consistent with the \textit{border rank} property. This property, peculiar to tensors of order three and higher, refers to the topological closure properties of the space of tensor decompositions. When this space is not closed, the abstract closure defines tensor decompositions on the ``border" of the space that have lower rank. In more concrete terms, it gives rise to approximate decompositions of lower rank, where the approximation as measured by $\epsilon$ gets better when the weights are allowed to diverge without limit.

\begin{figure}[t]
\begin{center}
\includegraphics[width=4.in]{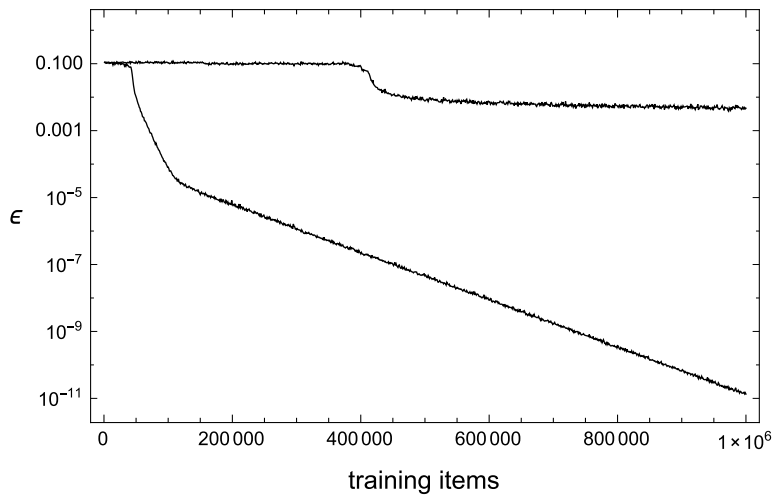}
\end{center}
\caption{Root-mean-square error $\epsilon$, in rank 23 decompositions of the $3\times 3$ matrix multiplication tensor, on two runs. The absence of convincing linear convergence in one of the runs is accompanied (Figure 4) by a slow growth in the maximum weight magnitude.}
\end{figure}

\begin{figure}[h!]
\begin{center}
\includegraphics[width=4.in]{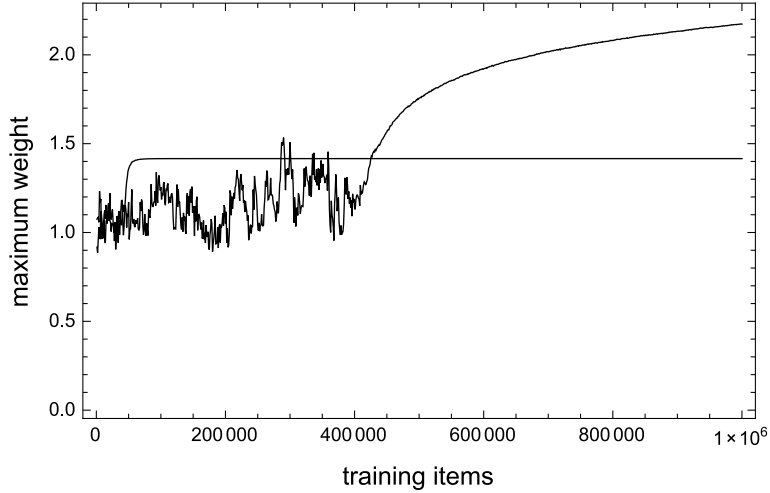}
\end{center}
\caption{The maximum weight magnitudes, $\max{(|W_a|,|W_b|,|W_c|)}$, for the two training runs shown in Figure 3.}
\end{figure}

We see evidence of ``border" behavior when we compare the evolution of the maximum weight magnitude, $\max{(|W_a|,|W_b|,|W_c|)}$, for the two runs in Figure 3. These are plotted in Figure 4. In contrast to the linearly convergent run, where the maximum weight saturates, the run with slow convergence has a corresponding slow growth in the maximum weight. An interpretation of our results, consistent with what is known about tensor decompositions of $M_3$, is that the weights in the slowly converging runs are rank 23 border approximations to true rank 24 or higher decompositions. Running the network with 22 multipliers we only observe the slow convergence/growing weight behavior. While this does not prove the true rank is 23, it is consistent with $M_3$ having rank 22 border approximations \cite{Sch}. That slow convergence was never observed in our experiments with $M_2$ is in agreement with the known fact \cite{JML} that this tensor does not have lower rank border approximations.

Conservative learning with the Strassen network for $3\times 3$ matrices is summarized in Figure 5. This gives the distribution of the final decomposition error $\epsilon$ for $10^3$ runs, each limited to $10^8$ training items and terminated when the error dropped below $10^{-14}$. Runs with linear convergence, such as the one in Figure 3, typically required far fewer than $10^8$ training items. The large peak at the low end of the distribution in Figure 5, about 64\% of all runs, therefore gives the probability that the network finds a true rank 23 decomposition. As explained above, we believe the network is slowly converging to border approximations in the other runs.

\begin{figure}[t]
\begin{center}
\includegraphics[width=4.in]{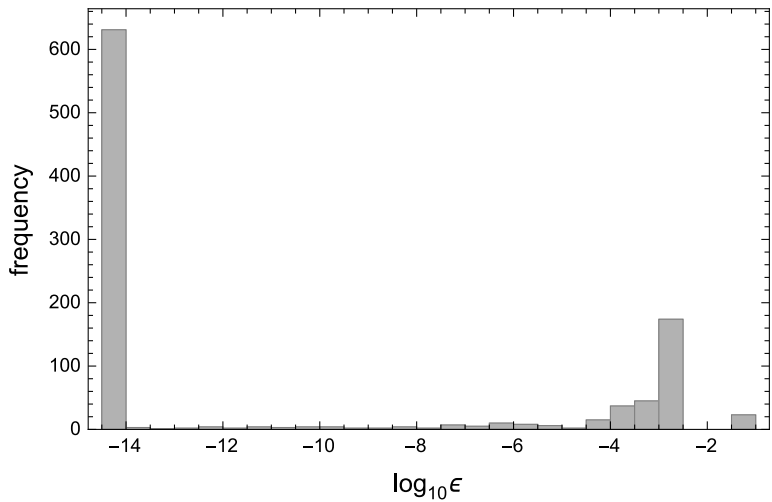}
\end{center}
\caption{Distribution of the final decomposition error of $M_3$ in $10^3$ runs, each with $10^8$ training items. Runs were terminated when $\epsilon$ dropped below $10^{-14}$.}
\end{figure}

\section{Appendix}

The starting point for deriving the conservative learning rules for the Strassen network (Fig. 1) is the Lagrangian function
\begin{multline}
\mathcal{L}=\frac{1}{2}\Tr{(\Delta_a^T\Delta_a)}+\frac{1}{2}\Tr{(\Delta_b^T\Delta_b)}+\frac{1}{2}\Tr{(\Delta_c^T\Delta_c)}\\
+\alpha^T\left(a^*-(W_a+\Delta_a)a\right)
+\beta^T\left(b^*-(W_b+\Delta_b)b\right)\\
+\gamma^T\left(c-(W_c+\Delta_c)a^*\circ b^*\right).
\end{multline}
To streamline the derivation we use matrix notation. The new training item is the triplet of vectors $(a,b,c)$ corresponding to the unrolled matrices $(A,B,C=AB)$. We should think of $(a,b)$ as inputs to the network and $c$ as the output. Prior to this item the network has weight matrices $W_a$ and $W_b$ that map the inputs to the pairs of inputs of the multipliers in the middle layer, $(a^*,b^*)$. The outputs of the multipliers is the vector $c^*=a^*\circ b^*$, where $\circ$ denotes componentwise multiplication. If the network weights are correct even for the new item, then mapping $c^*$ with the weight matrix $W_c$ should match the output vector $c$. In general, the three weight matrices have to be changed by $(\Delta_a,\Delta_b,\Delta_c)$ for this to be true, and the first three terms of $\mathcal{L}$ are the Frobenius norm objective on these changes to keep them small. The last three terms are constraints imposed via three vectors of Lagrange multipliers, $(\alpha,\beta,\gamma)$. These insure that with the changed weights the inputs/outputs of the network match the inputs/outputs of the multipliers in the middle layer.

Given the current weights $(W_a,W_b,W_c)$ and the new training item $(a,b,c)$, our task is then to find a stationary point of Lagrangian $\mathcal{L}$ for the variables $(\Delta_a,\Delta_b,\Delta_c)$, $(a^*,b^*)$, and $(\alpha,\beta,\gamma)$. In the derivation below we assume that the inputs are normalized as $a^T a=b^T b=1$.

Stationarity with respect to $\Delta_a$ and $\Delta_b$ imply
\begin{equation}\label{dadb}
\Delta_a=\alpha\, a^T\qquad \Delta_b=\beta\, b^T.
\end{equation}
Comparing with the update rule for the linear network in section \ref{sec:CL} prompts us to interpret $\alpha$ and $\beta$ as discrepancies associated with the multiplier inputs.
Multiplying \eqref{dadb} on the right respectively by $a$ and $b$ we also have
\begin{equation}\label{alphabeta}
\Delta_a\, a=\alpha\qquad \Delta_b\, b=\beta.
\end{equation}
These equations and ones to follow are consistent with the Lagrange multipliers vanishing proportionately with the changes in the weights.

Before we impose stationarity with respect to the other variables, we define a set of approximate vectors associated with implementing the network using the current (unchanged) weights. This mode of evaluating the nodes (vectors) in the network is called a forward pass.
\begin{equation}\label{fpass1}
\tilde{a}^*=W_a\, a\qquad \tilde{b}^*=W_b\, b.
\end{equation}
\begin{equation}\label{fpass2}
\tilde{c}=W_c\; \tilde{a}^*\circ \tilde{b}^*.
\end{equation}
Imposing stationarity with respect to $(\alpha,\beta)$ and comparing the resulting equations with \eqref{fpass1} and \eqref{alphabeta}, we obtain
\begin{equation}
a^*=(W_a+\Delta_a)a\qquad b^*=(W_b+\Delta_b)b,
\end{equation}
\begin{equation}
a^*-\tilde{a}^*=\Delta_a\,a=\alpha\qquad b^*-\tilde{b}^*=\Delta_b\,b=\beta.
\end{equation}
This shows that the discrepancies represented by $\alpha$ and $\beta$ are the differences between their true values (after the conservative-learning update) and their forward-pass values.

Learning in some sense starts at the output layer, where the $\tilde{c}$ of the forward pass is compared with the output $c$ of the training item. To obtain the conservative learning rule for this we start by imposing stationarity of $\mathcal{L}$ with respect to $\Delta_c$:
\begin{equation}\label{dcexact}
\Delta_c=\gamma\, (a^*\circ b^*)^T.
\end{equation}
We now make the first of a series of approximations. Extending to the output layer the property of the middle layer, that the pairs $(\Delta_a,\alpha)$ and $(\Delta_b,\beta)$ have the same order of smallness, we expect $(\Delta_c,\gamma)$ also to vanish proportionately. A good approximation of \eqref{dcexact} is then to replace $a^*$ and $b^*$ by their forward pass values (the error being higher order),
\begin{equation}\label{dc}
\Delta_c\approx\gamma\, (\tilde{a}^*\circ \tilde{b}^*)^T=\gamma\, \tilde{c}^{*T},
\end{equation}
where we have expressed the result in terms the forward pass value of the multiplier outputs, $\tilde{c}^*$.
Now imposing stationarity of $\mathcal{L}$ with respect to $\gamma$,
\begin{equation}
c=(W_c+\Delta_c)\,a^*\circ b^*,
\end{equation}
and comparing with \eqref{fpass2} we obtain
\begin{eqnarray}
c-\tilde{c}&\approx& \Delta_c\;\tilde{c}^*+W_c\left((a^*-\tilde{a}^*)\circ \tilde{b}^*+\tilde{a}^*\circ (b^*-\tilde{b}^*)\right)\\
&\approx& \gamma\, \tilde{c}^{*T} \tilde{c}^*+W_c\, (\alpha\circ\tilde{b}^*+\tilde{a}^*\circ \beta).\label{disc}
\end{eqnarray}
where we have neglected second order terms and used \eqref{dadb} and \eqref{dc}.

By imposing stationarity with respect to the two remaining sets of variables $(a^*,b^*)$ we will obtain equations that relate $(\alpha, \beta)$ to $\gamma$, enabling us to cast the equation for the output discrepancy \eqref{disc} just in terms of the unknown $\gamma$:
\begin{equation}
\alpha=b^*\circ (W_c^T+\Delta_c^T)\gamma\qquad \beta=a^*\circ (W_c^T+\Delta_c^T)\gamma,
\end{equation}
\begin{equation}\label{alphabetagamma}
\alpha=b^*\circ (W_c^T\gamma)\qquad \beta=a^*\circ (W_c^T\gamma).
\end{equation}
In \eqref{alphabetagamma} we again discarded second order terms. Substituting $\alpha$ and $\beta$ from \eqref{alphabetagamma} into \eqref{disc} we arrive at the equation that begins the process of updating the weights:
\begin{equation}\label{gammaeq}
c-\tilde{c}\approx  (\tilde{c}^{*T} \tilde{c}^*)\,\gamma+W_c\,(\tilde{a}^*\circ\tilde{a}^*+\tilde{b}^*\circ\tilde{b}^*)\,W_c^T\;\gamma.
\end{equation}

Equation \eqref{gammaeq} relates the discrepancy, between the true (training) $c$ and the $\tilde{c}$ of the forward pass, to the Lagrange multipliers $\gamma$. Were we to neglect the off-diagonal terms in this linear matrix equation and determine $\gamma$ by
\begin{equation}\label{simplebackprop}
\gamma \leftarrow \frac{1}{\tilde{c}^{*T} \tilde{c}^*}(c-\tilde{c}),
\end{equation}
then the process of learning item $(a,b,c)$ would be completely analogous to the usual ``back-propagation" scheme. After forward-propagating the inputs $(a,b)$ with $W_a$, $W_b$ and $W_c$ to determine the discrepancy $c-\tilde{c}$, the $\gamma$ given by \eqref{simplebackprop} is back-propagated \eqref{alphabetagamma} by $W_c^T$ to get $\alpha$ and $\beta$. The three Lagrange multipliers then determine the weight changes by \eqref{dadb} and \eqref{dc}.

Lacking an argument for discarding the off-diagonal terms in \eqref{gammaeq}, we need to look for methods to solve this more complex linear equation. The off-diagonal terms correspond to one-level of back-propagation followed by forward-propa\-gation. An iterative solution of the linear equation would thus involve multiple backward-forward propagations between just the final two layers of the network. This is not as intimidating as it might seem for two reasons. First, when using the conjugate gradient (CG) method, the number of backward-forward iterations is bounded by the number of components of the vector $\gamma$. Second, given a reasonable initial solution-estimate, CG usually requires only few iterations in practice.

Since our derivation of the conservative learning rules has been based on the premise that the discrepancies $(a^*-\tilde{a}^*, b^*-\tilde{b}^*, c-\tilde{c})$ are small, we are keeping with this premise when we take as our initial CG solution-estimate $\gamma=0$. As a better motivated alternative to \eqref{simplebackprop} we apply the fewest (non-trivial) number of CG iterations --- a single one --- to this solution-estimate. At this single-iteration level of CG the approximate solution has a simple interpretation. Let
\begin{equation}\label{G}
G=(\tilde{c}^{*T} \tilde{c}^*)+W_c\,(\tilde{a}^*\circ\tilde{a}^*+\tilde{b}^*\circ\tilde{b}^*)\,W_c^T
\end{equation}
be the symmetric positive definite matrix in our linear equation
\begin{equation}
\delta=c-\tilde{c}=G\gamma.
\end{equation}
The single-iteration CG solution has the form $\gamma_1=\lambda \delta$, where $\lambda$ is a scalar multiplier and determined by projecting the equation on $\delta$:
\begin{equation}
\delta^T(\delta-G(\lambda\delta))=0.
\end{equation}
The conservative learning alternative to \eqref{simplebackprop} is therefore
\begin{equation}\label{CGprop}
\gamma \leftarrow \gamma_1=\frac{(c-\tilde{c})^T(c-\tilde{c})}{(c-\tilde{c})^T G (c-\tilde{c})}(c-\tilde{c}).
\end{equation}
Not surprisingly this reduces to \eqref{simplebackprop} when the off-diagonal terms in $G$ are dropped.

We conclude this appendix with a short summary of the conservative learning rules for updating the network weights when given a new training item $(a,b,c)$.
\begin{enumerate}
\item Forward propagate $(a,b)$ using \eqref{fpass1} and \eqref{fpass2} to get $(\tilde{a}^*,\tilde{b}^*, \tilde{c}^*= \tilde{a}^*\circ  \tilde{b}^*)$ and $\tilde{c}$.
\item From the output discrepancy $c-\tilde{c}$ compute $\gamma$ using \eqref{CGprop}. Computing the scalar multiple of $c-\tilde{c}$ in this expression requires a single backward-forward propagation by the matrix $G$ in \eqref{G}.
\item The update $\Delta_c$ of $W_c$ is the rank 1 matrix \eqref{dc} constructed from $\gamma$ and $\tilde{c}^*$.
\item Backward propagate $\gamma$ by \eqref{alphabetagamma} to obtain $\alpha$ and $\beta$.
\item The updates $(\Delta_a,\Delta_b)$ of $(W_a,W_b)$ are given by the rank 1 matrices \eqref{dadb} constructed from $(\alpha,\beta)$ and $(a,b)$.
\end{enumerate}

\section{Acknowledgements}\label{sec:ack}

I thank Cris Moore for instigating this study, Jonathan Yedidia for suggesting the conservative learning method, and Alex Alemi for keeping things competitive. The Simons Foundation and ADI Lyric Labs provided financial support.

\end{document}